\documentclass[fleqn,a4paper,12pt]{article}
\usepackage{amsmath,amsfonts,calrsfs,amssymb,color}

\newtheorem{Theorem}{Theorem}[section]

\newtheorem{Proposition}[Theorem]{Proposition}

\newtheorem{Corollary}[Theorem]{Corollary}
\newtheorem{Remark}[Theorem]{Remark}
\newtheorem{Example}[Theorem]{Example}

\makeatletter
\@addtoreset{equation}{section}

\makeatother

\def\R{\mathbb R}
\def\N{\mathbb N}
\def\Z{\mathbb Z}
\def\C{\mathbb C}
\def\E{\mathbb E}
\def\P{\mathbb P}
\def\ds{\displaystyle}
 \def\oo{\mathaccent23}

\newcommand{\one}{1\!\!\!\;\mathrm{l}}

\title{\bf Singular Stochastic Equations on Hilbert Spaces: Harnack Inequalities for their Transition Semigroups}
\author{
Giuseppe Da Prato \thanks{Supported in part by ``Equazioni di
Kolmogorov'' from the Italian
``Ministero della Ricerca Scientifica e Tecnologica''
} ,\\
 Scuola Normale Superiore
di Pisa, Italy\\
 \\
Michael R\"ockner \thanks{Supported by the DFG through SFB-701 and IRTG 1132, by NSF-Grant 0603742 as well as by the
BIBOS-Research Center. } \\
Faculty of Mathematics, University of Bielefeld, Germany\\ and\\
Department of Mathematics and Statistics,\\ Purdue University, W. Lafayette, 47906, IN,   U. S. A.\\\\
Feng-Yu Wang \thanks{The corresponding author. Supported in part by  WIMCS, Creative Research Group Fund of the National Natural Science Foundation of China (No. 10721091) and the 973-Project.}\\
 School of Math. Sci.  and  Lab. Math Com. Sys.\\ Beijing Normal University,  100875, China\\
and
\\
  Department of Mathematics, Swansea University,\\ Singleton Park, SA2 8PP, Swansea, UK
}
\date{ }
\begin{document}
\maketitle

 \begin{abstract}
We consider stochastic equations in Hilbert spaces with singular
drift in the framework of \cite{4}. We prove a Harnack inequality
(in the sense of \cite{W97}) for its transition semigroup and
exploit its consequences. In particular, we prove regularizing and ultraboundedness properties of the transition semigroup as well as  that  the corresponding Kolmogorov operator has at most one infinitesimally invariant measure $\mu$ (satisfying some mild integrability conditions). Finally, we prove existence of such a measure $\mu$ for non-continuous drifts.
\end{abstract}
\bigskip

\noindent {\bf 2000 Mathematics Subject Classification AMS}: 60H15, 35R15,  35J15,

 \bigskip

\def\R{\mathbb R}  \def\ff{\frac} \def\ss{\sqrt} \def\BB{\mathbb
B}
\def\N{\mathbb N} \def\kk{\kappa} \def\m{{\bf m}}
\def\dd{\delta} \def\DD{\Delta} \def\vv{\varepsilon} \def\rr{\rho}
\def\<{\langle} \def\>{\rangle} \def\GG{\Gamma} \def\gg{\gamma}
  \def\nn{\nabla} \def\pp{\partial} \def\tt{\tilde}
\def\d{\text{\rm{d}}} \def\bb{\beta} \def\aa{\alpha} \def\D{\scr D}
\def\E{\mathbb E} \def\si{\sigma} \def\ess{\text{\rm{ess}}}
\def\beg{\begin} \def\beq{\begin{equation}}  \def\F{\scr F}
\def\Ric{\text{\rm{Ric}}} \def\Hess{\text{\rm{Hess}}}\def\B{\mathbb B}
\def\e{\text{\rm{e}}} \def\ua{\underline a} \def\OO{\Omega} \def\sE{\scr E}
\def\oo{\omega}     \def\tt{\tilde} \def\Ric{\text{\rm{Ric}}}
\def\cut{\text{\rm{cut}}} \def\P{\mathbb P} \def\ifn{I_n(f^{\bigotimes n})}
\def\C{\scr C}      \def\aaa{\mathbf{r}}     \def\r{r}
\def\gap{\text{\rm{gap}}} \def\prr{\pi_{{\bf m},\varrho}}  \def\r{\mathbf r}
\def\Z{\mathbb Z} \def\vrr{\varrho} \def\ll{\lambda}
\def\L{\scr L}\def\Tt{\tt} \def\TT{\tt} \def\C{\mathbb C}

\section{Introduction, framework and main results}

In this paper we continue our study of stochastic equations in
Hilbert spaces with singular drift through its associated Kolmogorov
equations started in \cite{4}. The main aim is to prove a Harnack
inequality for its transition semigroup in the sense of \cite{W97}
(see also \cite{1,10,12} for further development) and exploit its
consequences. See also \cite{OY} for an improvement of the main results in \cite{10}
concerning  generalized Mehler semigroups. To describe our results more
 precisely, let us first
recall the framework from \cite{4}.

Consider the stochastic equation
\begin{equation}
\label{e1.1}
\left\{\begin{array}{l}
dX(t)=(AX(t)+F(X(t)))dt+\sigma dW(t)\\
\\
X(0)=x\in H.
\end{array}\right.
\end{equation}
Here $H$ is a real separable Hilbert space with inner product $\langle\cdot ,\cdot  \rangle$ and norm $|\cdot|$, $W=W(t),\;t\ge 0, $ is a cylindrical Brownian motion on $H$ defined on a stochastic basis $(\Omega, \mathcal F, (\mathcal F_t)_{t\ge 0}, \P)$ and the coefficients satisfy the following hypotheses:\bigskip

\noindent (H1) $(A,D(A))$ is the generator of a $C_0$-semigroup, $T_t=e^{tA}$, $t\ge 0$, on $H$ and for some $\omega\in \R$
\begin{equation}
\label{e1.1'}
\langle Ax,x  \rangle\le\omega |x|^2,\quad\forall\;x\in D(A).
 \end{equation}
\bigskip

\noindent (H2) $\sigma \in L(H)$ (the space of all bounded linear  operators on $H$) such that $\sigma $ is positive definite, self-adjoint and
\begin{enumerate}
\item[(i)] $\ds\int_0^\infty (1+t^{-\aa})\|T_t\sigma\|^2_{HS}dt<\infty$ for some $\alpha >0,$
where $\|\cdot \|_{HS}$ denotes the norm on the space of all
Hilbert--Schmidt operators on $H$.

\item[(ii)] $\sigma ^{-1}\in L(H)$. \bigskip
\end{enumerate}

\noindent (H3)  $F:D(F)\subset H\to 2^H$ is an $m$-dissipative map, i.e.,
$$
\langle u-v,x-y  \rangle\le 0,\quad\forall\;x,y\in D(F),\;u\in F(x),\;v\in F(y),
$$
(``dissipativity'') and
$$
\mbox{\rm Range}\;(I-F):=\bigcup_{x\in D(F)} (x-F(x))=H.
$$
Furthermore, $F_0(x)\in F(x),\;x\in D(F),$ is such that
$$
|F_0(x)|=\min_{y\in F(x)}|y|.
$$

Here we recall that for $F$ as in (H3) we have that $F(x)$ is closed, non empty and convex.

The corresponding Kolmogorov operator is then given as follows: Let $\mathcal E_A(H)$
denote the linear span of all real parts  of functions of the form $\varphi=e^{i\langle h,\cdot
 \rangle}$, $h\in D(A^*)$, where $A^*$ denotes the adjoint operator of $A$, and define for any $x\in D(F)$,
\begin{equation}
\label{e1.2}
L_0\varphi(x)=\frac12\;\mbox{\rm Tr}\;(\sigma^2D^2\varphi(x))+\langle x, A^*D\varphi(x)   \rangle+
\langle F_0(x),D\varphi(x)   \rangle,\quad \varphi\in \mathcal E_A(H).
\end{equation}
Additionally, we assume:\bigskip

\noindent (H4) There exists a probability measure $\mu$ on $H$ (equipped with its Borel
$\sigma$-algebra $\mathcal B(H)$) such that
\begin{enumerate}
\item[(i)] $\mu(D(F))=1$,

\item[(ii)] $\ds\int_H (1+|x|^2)(1+|F_0(x)|)\mu(dx)<\infty,$

\item[(iii)] $\ds \int_H L_0\varphi d\mu=0$ for all
$\varphi\in \mathcal E_A(H)$.
\end{enumerate}

\begin{Remark}
\label{r1.1}
\em (i) A measure for which the last equality in (H4) (makes sense and) holds is called
{\em infinitesimally invariant} for $(L_0,\mathcal E_A(H))$.\bigskip

\noindent (ii) Since $\omega$ in \eqref{e1.1'} is an arbitrary real number we can relax (H3) by allowing that for some $c\in (0,\infty)$
$$
\langle u-v,x-y\rangle\le c |x-y|^2,\quad\forall\;x,y\in D(F),\;u\in F(x),\;v\in F(y).
$$
We simply replace $F$  by $F-c$ and $A$ by $A+c$ to reduce this case to (H3).\bigskip

\noindent (iii) At this point we would like to stress that under the above assumptions (H1)-(H4) (and (H5) below)
because $F_0$ is merely measurable and $\sigma$ is not   Hilbert-Schmidt, it is unknown whether \eqref{e1.1}
has a strong solution. \bigskip

\noindent (iv) Similarly as in \cite{4} (see \cite[Remark  4.4]{4} in particular) we expect that (H2)(ii) can be relaxed to the  condition that $\sigma=(-A)^{-\gamma}$ for some $\gamma\in [0,1/2]$. However, some of the approximation arguments below become more involved. So, for simplicity we assume (H2)(ii).
\end{Remark}

The following are the main results of \cite{4} which we shall use below.
\begin{Theorem}
\label{t1.2} {\bf (cf. [6, Theorem 2.3 and Corollary
2.5])} Assume $(H1)$, $(H2)(i)$, $(H3)$ and $(H4)$.
Then for any measure $\mu$ as in $(H4)$ the operator $(L_0,\mathcal
E_A(H))$ is dissipative on $L^1(H,\mu)$, hence closable. Its closure
$(L_\mu,D(L_\mu))$ generates a $C_0$-semigroup $P_t^\mu$, $t\ge 0$,
on $L^1(H,\mu)$ which is Markovian, i.e., $P_t^\mu 1=1$ and $P_t^\mu
f\ge 0$ for all nonnegative $f\in L^1(H,\mu)$ and all $t>0$.
Furthermore, $\mu$ is $P_t^\mu$-invariant, i.e.,
$$
\int_HP_t^\mu fd\mu=\int_H fd\mu,\quad\forall\;f\in
L^1(H,\mu).
$$
\end{Theorem}
Below $B_b(H)$, $C_b(H)$ denote the bounded Borel-measurable, continuous functions respectively from $H$ into
$\R$ and $\|\cdot\|$ denotes the usual norm
on $L(H)$.
\begin{Theorem} {\bf (cf. [6, Proposition  5.7])}
\label{t1.3}
  Assume $(H1)$--$(H4)$ hold. Then for any measure $\mu$ as in
 $(H4)$ and $H_0:=\mbox{\rm supp}\;\mu$ (:=largest closed set of $H$ whose complement is a
 $\mu$-zero set) there exists  a semigroup $p_t^\mu(x,dy)$, $x\in H_0$, $t>0$, of kernels such that
 $p_t^\mu f$ is a $\mu$-version of $P_t^\mu f$ for all $f\in B_b(H)$, $t>0,$ where as usual
$$
p_t^\mu f(x)=\int_H f(y)p_t^\mu(x,dy),\quad x\in H_0.
$$
Furthermore, for all $f\in B_b(H)$, $t>0,$ $x,y\in H_0$
\begin{equation}
\label{e1.3}
|p_t^\mu f(x)-p_t^\mu f(y)|\le \frac{e^{|\omega|t}}{
\sqrt{t\wedge 1}}\; \|f\|_0 \|\sigma^{-1}\||x-y|
\end{equation}
and for all $f\in Lip_b(H)$ (:= all bounded Lipschitz functions on $H$)
\begin{equation}
\label{e1.4}
|p_t^\mu f(x)-p_t^\mu f(y)|\le e^{|\omega|t} \|f\|_{Lip} |x-y|,\quad\forall\;t>0,\; x,y\in H_0,
\end{equation}
and
\begin{equation}
\label{e1.5}
\lim_{t\to 0}p_t^\mu f(x)=f(x),\quad\forall\; x\in H_0.
\end{equation}
$($Here $\|f\|_0 $, $\|f\|_{Lip}$ denote the supremum, Lipschitz norm of $f$ respectively.$)$
 Finally, $\mu$ is $p_t^\mu$-invariant.
\end{Theorem}
\begin{Remark}
\label{r1.4}
\em (i) Both results above have been proved in \cite{4} on $L^2(H,\mu)$ rather than on $L^1(H,\mu)$,
but the proofs for $L^1(H,\mu)$ are entirely analogous.

(ii) In \cite{4} we assume $\omega$ in (H1) to be negative, getting a stronger estimate than \eqref{e1.3}
(cf. \cite[(5.11)]{4}). But the same proof as in \cite {4} leads to \eqref{e1.3} for arbitrary
$\omega\in \R$ (cf. the proof of \cite[Proposition 4.3]{4} for $t\in [0,1]$). Then by virtue of the
semigroup property and since $p_t^\mu$ is Markov we get \eqref{e1.3} for all $t>0$.

(iii) Theorem \ref{t1.3} holds in more general situations since (H2)(ii) can be relaxed (cf. \cite[Remark 4.4]{4}
and \cite[Proposition 8.3.3]{3a}).

(iv) \eqref{e1.3} above implies that $p_t^\mu$, $t>0,$ is strongly
Feller, i.e., $p_t^\mu(B_b(H))\subset C(H_0)$ (=all continuous
functions on $H_0$). We shall prove below that  under the additional
condition (H5) we even have $p_t^\mu(L^p(H,\mu))\subset C(H_0)$ for
all $p>1$ and that $\mu$ in (H4) is unique. However, so far we have
not been able to prove  that for this unique $\mu$ we have supp
$\mu=H$, though we conjecture that this is true. $\Box$

\end{Remark}\bigskip

For the results on Harnack inequalities, in this paper we need one more condition.\bigskip

  \begin{enumerate}

   \item[(H5) (i)]   $(1+\omega-A,D(A))$ satisfies the weak sector condition (cf. e.g. \cite{7a}), i.e.,
   there exists a constant $K>0$ such that
\begin{equation}
\label{e1.6}
\langle (1+\omega-A)x,y  \rangle\le K
\langle (1+\omega-A)x,x  \rangle^{1/2}\langle (1+\omega-A)y,y  \rangle^{1/2},\quad\forall\;x,y\in D(A).
\end{equation}
\item[(ii)] There exists a sequence of $A$-invariant
finite dimensional subspaces $H_n\subset D(A)$  such that $\bigcup_{n=1}^\infty H_n$ is dense in $H$.
\end{enumerate}
We note that if $A$ is self-adjoint, then (H2) implies that $A$ has a discrete spectrum   which in turn implies
that (H5)(ii) holds.
\begin{Remark}
\label{r1.5}
\em Let $(A,D(A))$ satisfy (H1). Then the following is well known:

(i) (H5)  (i) is equivalent to the fact that the semigroup generated by $(1+\omega-A,D(A))$ on the complexification
$H_\C$ of $H$ is a holomorphic contraction semigroup on $H_\C$ (cf. e.g. \cite[Chapter I, Corollary 2.21]{7a}).

(ii) (H5) (i) is equivalent to $(1+\omega-A,D(A))$ being variational. Indeed, let $(\mathcal E, D(\mathcal E))$ be
the coercive closed form generated by $(1+\omega-A,D(A))$ (cf.  \cite[Chapter I, Section 2]{7a})
and $(\widetilde{\mathcal E}, D(\mathcal E))$ be
its symmetric part. Then define
\begin{equation}
\label{e1.7}
V:=D(\mathcal E)\;\mbox{\rm with inner product}\;\widetilde{\mathcal E}\;\mbox{\rm and $V^*$ to be its dual}.
\end{equation}
Then
\begin{equation}
\label{e1.8}
V\subset H\subset V^*
\end{equation}
and $1+\omega-A: D(A)\to H$ has a natural unique continuous extension from $V$ to $V^*$ satisfying all the required
properties  (cf.  \cite[Chapter I, Section 2, in particular Remark 2.5]{7a}).

\end{Remark}
Now we can formulate the main result of this paper, namely the Harnack inequality for $p_t^\mu,\;t>0.$

\begin{Theorem}
\label{t1.6}
Suppose $(H1)-(H5)$ hold and let $\mu$ be any measure as in $(H4)$  and $p_t^\mu(x,dy)$ as in Theorem $\ref{t1.3}$ above.
Let $p\in (1,\infty)$.
Then for all $f\in B_b(H)$, $f\ge 0,$
\begin{equation}
\label{e1.9} (p_t^\mu f(x))^p\le p_t^\mu f^p(y)
\exp\left[\|\sigma^{-1}\|^2\;\frac{p\omega|x-y|^2}{(p-1)
(1-e^{-2\omega t})}
 \right],\quad t>0,\;x,y\in H_0.
\end{equation}

\end{Theorem}
As   consequences in the situation of Theorem \ref{t1.6} (i.e. assuming (H1)-(H5)) we obtain:
\begin{Corollary}
\label{c1.7} For all $t>0$ and $p\in (1,\infty)$
$$
p^\mu_t(L^p(H,\mu))\subset C(H_0).
$$
\end{Corollary}
\begin{Corollary}
\label{c1.8}
$\mu$ in $(H4)$ is unique.
\end{Corollary}
Because of this result below we write $p_t(x,dy)$ instead of
$p_t^\mu(x,dy)$.

Finally, we have
\begin{Corollary}
\label{c1.9} (i) For every $x\in H_0$, $p_t(x,dy)$ has a density
$\rho_t(x,y)$ with respect to $\mu$ and
\begin{equation}
\label{e1.10}
\|\rho_t(x,\cdot)\|^{p/(p-1)}_p\le \frac{1}{ \int_H\exp\left[-\|\sigma^{-1}\|^2\;\frac{p\omega|x-y|^2}
{ (1-e^{-2\omega t})}\right]\mu(dy)   },\quad  x\in H_0,\; p\in (1,\infty).
\end{equation}
(ii) If $\mu(e^{\lambda|\cdot|^2})<\infty$ for some
$\lambda>2(\omega\wedge 0)^2\|\sigma^{-1}\|^2,$ then $p_t$ is
hyperbounded, i.e. $\|p_t\|_{L^2(H,\mu)\to L^4(H,\mu)}<\infty$ for
some $t>0$.
\end{Corollary}

\begin{Corollary}
\label{c1.10}
For simplicity, let $\sigma=I$ and
 instead of $(H1)$ assume that more strongly $(A,D(A))$ is self-adjoint satisfying
 \eqref{e1.1'}. We furthermore assume that $|F_0|\in L^2(H,\mu).$

 \noindent (i) There exists $M\in \mathcal B(H_0)$, $M\subset D(F)$, $\mu(M)=1$ such that for every $x\in M$ equation \eqref{e1.1} has a pointwise unique continuous strong solution (in the mild sense see \eqref{e4.8ro} below), such that $X(t)\in M$ for all $t\ge 0$ $\P$-a.s..

\noindent (ii) Suppose there exists $\Phi\in C([0,\infty))$
positive and strictly increasing such that
$\lim_{s\to\infty}s^{-1}\Phi(s)=\infty$ and
\begin{equation}
\label{L1}
 \Psi(s):= \int_s^\infty \ff{\d r}{\Phi(r)}<\infty,\quad\forall\; s>0.
\end{equation}
 If there exists a constant $c>0$ such that
\begin{equation}
\label{L2} \langle F_0(x)-F_0(y), x-y\rangle\le
c-\Phi(|x-y|^2),\quad \forall\; x,y\in D(F),
\end{equation}
 then $p_t$ is ultrabounded with
$$\|p_t\|_{L^2(H,\mu)\to L^\infty(H,\mu)}\le \exp\Big[\ff{\ll (1+
\Psi^{-1}(t/4))}{(1-\e^{-\oo t/2})^2}\Big],\ \ \ t>0,$$
 holding for
some constant $\ll>0.$
 \end{Corollary}

\begin{Remark}
\label{r1.10}
 \em

 We emphasize that since the nonlinear part
$F_0$ of our Kolmogorov operator is in general not continuous, it
was quite surprising for us that in this infinite dimensional case
nevertheless the generated semigroup $P_t$ maps $L^1$- functions to
continuous ones as stated in Corollary  \ref{c1.7}.

\end{Remark}

The proof that Corollary \ref{c1.9} follows from Theorem \ref{t1.6} is completely standard. So, we will omit the proofs
and instead refer to \cite{10}, \cite{12}.

Corollary \ref{c1.7} is new and follows whenever a semigroup $p_t$ satisfies the Harnack inequality
(see Proposition \ref{p4.1} below).

Corollary \ref{c1.8} is   new. Since \eqref{e1.9} implies irreducibility of $p^\mu_t$ and Corollary \ref{c1.7}
implies that it is strongly Feller, a well known theorem due to Doob immediately implies that $\mu$ is the unique
invariant measure for $p^\mu_t,\;t>0.$ $p^\mu_t$, however, depends on $\mu$,
so Corollary \ref{c1.8} is a  stronger statement.
Corollary \ref{c1.10} is also new.

Theorem \ref{t1.6} as well as Corollaries \ref{c1.7}, \ref{c1.8} and
\ref{c1.10} will be proved in     Section 4. In Section 3 we first prove Theorem \ref{t1.6} in case $F_0$ is Lipschitz, and in Section 2 we prepare the tools that allow us to reduce the general case to the Lipschitz case.
In Section
5 we prove two results (see Theorems \ref{t5.2} and \ref{t5.4}) on the existence of a measure satisfying (H4) under some
additional conditions and present an application to an example where $F_0$ is not continuous.
For a discussion of a number of other explicit examples  satisfying our conditions see \cite[Section 9]{4}.\bigskip

\section{Reduction to regular $F_0$}

Let $F$ be as in (H3). As in \cite{4} we may consider the Yosida approximation of $F$, i.e., for any $\alpha>0$ we
set
\begin{equation}
\label{e2.1}
F_\alpha(x):=\frac1\alpha\;(J_\alpha(x)-x),\quad x\in H,
\end{equation}
where for $x\in H$
$$
J_\alpha(x):=(I-\alpha F)^{-1}(x),\quad  \alpha>0,
$$
and $I(x):=x$. Then each $F_\alpha$ is single valued, dissipative and it is well known that
\begin{equation}
\label{e2.2}
\lim_{\alpha\to 0}F_\alpha(x)=F_0(x),\quad\forall\;x\in D(F),
\end{equation}
\begin{equation}
\label{e2.3}
|F_\alpha(x)|\le |F_0(x)|,\quad\forall\;x\in D(F).
\end{equation}
Moreover, $F_\alpha$ is Lipschitz continuous, so $F_0$ is $\mathcal B(H)$-measurable. Since $F_\alpha$ is not
differentiable in general, as in \cite{4} we introduce a further regularization by setting
\begin{equation}
\label{e2.4}
F_{\alpha,\beta}(x):=\int_He^{\beta B}F_\alpha(e^{\beta B}x+y)N_{\frac12\;B^{-1}(e^{2\beta B}-1)}(dy),
\quad\alpha,\beta>0,
\end{equation}
where $B:D(B)\subset H\to H$ is a self-adjoint, negative definite linear operator such that $B^{-1}$
is of trace class and as usual for a trace class operator $Q$ the measure $N_Q$ is just the standard centered
Gaussian measure with covariance given by $Q$.

$F_{\alpha,\beta}$ is dissipative, of class $C^\infty$, has bounded derivatives of all the orders and
$F_{\alpha,\beta}\to F_{\alpha}$ pointwise as $\beta\to 0$.

Furthermore, for $\alpha>0$
\begin{equation}
\label{e2.5}
c_\alpha:=\sup\left\{\frac{|F_{\alpha,\beta}(x)|}{1+|x|}:\;x\in H,\;\beta\in (0,1]  \right\}<\infty.
\end{equation}
We refer to \cite[Theorem 9.19]{5} for details.

Now we consider the following regularized stochastic equation
\begin{equation}
\label{e2.6}
\left\{\begin{array}{l}
dX_{\alpha,\beta}(t)=(AX_{\alpha,\beta}(t)+F_{\alpha,\beta}(X_{\alpha,\beta}(t)))dt+\sigma dW(t)\\
\\
X_{\alpha,\beta}(0)=x\in H.
\end{array}\right.
\end{equation}
It is well known that \eqref{e2.6} has a unique mild solution $X_{\alpha,\beta}(t,x),\;t\ge 0$. Its associated
transition semigroup is given by
$$
P^{\alpha,\beta}_tf(x)=\E[f(X_{\alpha,\beta}(t,x))],\quad t>0,\;x\in H,
$$
for any $f\in B_b(H)$.  Here $\E$ denotes expectation with respect to $\P$.
\begin{Proposition}
\label{p2.1}
Assume $(H1)-(H4)$. Then there exists a $K_\sigma$-set $K\subset H$ such that $\mu(K)=1$ and for all $f\in B_b(H)$,
 $T>0$ there exist subsequences $(\alpha_n),\;(\beta_n)\to 0$ such that for all $x\in K$
\begin{equation}
\label{e2.7} \lim_{n\to\infty}\lim_{m\to\infty}
P^{\alpha_n,\beta_m}_\bullet f(x)= p^\mu_\bullet f(x)\quad \mbox{\it
weakly in}\;L^2(0,T;dt).
\end{equation}
\end{Proposition}
{\bf Proof}. This follows immediately from the proof of \cite[Proposition 5.7]{4}. (A closer look at the proof
even shows that \eqref{e2.7} holds for all $x\in H_0=$ supp $\mu$.) $\Box$

As we shall see in   Section 4, the proof of Theorem \ref{t1.6}   follows from Proposition \ref{p2.1} if we can prove
the corresponding Harnack inequality for each $P^{\alpha,\beta}_t$. Hence in the next section we confine ourselves to
the case when $F_0$ is dissipative and Lipschitz.

\section{The Lipschitz case}

In this section we assume that (H1)-(H3) and (H5)   hold and that
$F_0$   in (H3) is    in addition  Lipschitz continuous. The aim
of this section is to prove Theorem \ref{t1.6} for such special $F_0$ (see Proposition \ref{p3.1} below).
We shall do this by finite dimensional (Galerkin) approximations, since for the approximating finite dimensional
 processes we can apply the usual coupling argument.

We first note that since $F_0$ is Lipschitz \eqref{e1.1} has a unique mild solution $X(t,x),\;t\ge 0,$
for every initial condition $x\in H$ (cf.\cite{5})
and we denote the corresponding transition semigroup by $P_t$, $t>0$, i.e.
$$
P_tf(x):=\E[f(X(t,x))],\quad t>0,\; x\in X,
$$
where $f\in B_b(H)$.

Now we need to consider  an appropriate Galerkin approximation. To this end let $e_k\in D(A),\; k\in \N$, be
orthonormal such that $H_n=\mbox{\rm linear span}\;\{e_1,...,$ $e_n\},\; n\in \N.$ Hence
$\{ e_k:\;k\in \N\}$ is an orthonormal basis of $(H,\langle\cdot,\cdot   \rangle)$. Let $\pi_n:H\to H_n$ be
the orthogonal projection with respect to $(H,\langle\cdot,\cdot   \rangle)$. So, we can define
\begin{equation}
\label{e3.1}
A_n:=\pi_nA_{| H_n}\;(=A_{| H_n}\;\mbox{\rm by}\;(H5)(ii))
\end{equation}
and, furthermore
$$
F_n:=\pi_n F_{0| H_n},\quad\sigma_n:=\pi_n \sigma_{| H_n}.
$$
Obviously, $\sigma_n:H_n\to H_n$ is a self-adjoint, positive definite linear operator on $H_n$. Furthermore,
$\sigma_n$ is bijective, since it is one-to-one. To see the latter, one simply picks an orthonormal basis
$\{e_1^\sigma,...,e_n^\sigma\}$ of $H_n$ with respect to the inner product $\langle \cdot, \cdot  \rangle_\sigma$
defined by $\langle x, y  \rangle_\sigma:=\langle \sigma x, y  \rangle.$ Then if $x\in H_n$ is such that
$\sigma _nx=\pi_n\sigma  x=0$, it follows that
$$
\langle x, e_i^\sigma  \rangle_\sigma=
\langle \sigma x, e_i^\sigma  \rangle
=0,\quad\forall\;1\le i\le n.
$$
But $x=\sum_{i=1}^n\langle x, e_i^\sigma  \rangle_\sigma e_i^\sigma$, hence $x=0$.

Now fix $n\in \N$ and on $H_n$ consider the stochastic equation
\begin{equation}
\label{e3.2}
\left\{\begin{array}{l}
dX_n(t)=(A_nX_n(t)+F_n(X_n(t)))dt+\sigma_n dW_n(t)\\
\\
X_n(0)=x\in H_n,
\end{array}\right.
\end{equation}
where $W_n(t)=\pi_nW(t)=\sum_{i=1}^n\langle e_k, W(t) \rangle e_k.$

\eqref{e3.2} has a unique strong solution $X_n(t,x),\; t\ge 0$, for every initial condition $x\in H_n$ which
is pathwise continuous $\P$-a.s..
Consider the associated transition semigroup defined as before by
\begin{equation}
\label{e3.3}
P^n_tf(x)=\E[f(X_n(t,x))],\quad t>0,\;x\in H_n,
\end{equation}
where $f\in B_b(H_n)$.

Below we shall prove the following:
\begin{Proposition}
\label{p3.1}
Assume that $(H1)-(H5)$ hold. Then:

\noindent (i) For all $f\in C_b(H)$ and all $t>0$
$$
\lim_{n\to \infty}P^n_tf(x)=P_tf(x),\quad \forall x\in H_{n_0},\;n_0\in \N.
$$
\noindent (ii) For all nonnegative $f\in B_b(H)$ and all $n\in N$, $p\in (1,\infty)$
\begin{equation}
\label{e3.4}
(P_t^n f(x))^p\le P_t^nf^p(y) \exp\left[\|\sigma^{-1}\|^2\;\frac{p\omega|x-y|^2}{(p-1) (1-e^{-2\omega t})}
\right],\quad t>0,\;x,y\in H_n.
\end{equation}
\end{Proposition}
{\bf Proof}. (i): Define
$$
W_{A,\sigma}(t):=\int_0^te^{(t-s)A}\sigma dW(s),\quad t\ge 0.
$$
Note that by (H2)(i) we have that $W_{A,\sigma}(t),\;t\ge 0,$ is well defined and pathwise continuous.
 For $x\in H_{n_0},\;n_0\in \N$ fixed, let $Z(t),\;t\ge 0,$ be the unique variational solution
 (with triple $V\subset H\subset V^*$ as in Remark \ref{r1.5}(ii), see e.g. \cite{11a}) to
\begin{equation}
\label{e3.5}
\left\{\begin{array}{l}
dZ(t)=[AZ(t)+F_0(Z(t)+W_{A,\sigma}(t))]dt\\
Z(0)=x,
\end{array}\right.
\end{equation}
which then automatically satisfies
\begin{equation}
\label{e3.6}
\E\sup_{t\in[0,T]}|Z(t)|^2<+\infty.
\end{equation}
Then we have (see \cite{5}) that $Z(t)+W_{A,\sigma}(t),\;t\ge 0,$ is a mild solution to \eqref{e1.1}
(with $F_0$ Lipschitz), hence by uniqueness
\begin{equation}
\label{e3.7}
X(t,x)=Z(t)+W_{A,\sigma}(t),\quad t\ge 0.
\end{equation}
Clearly, since
\begin{equation}
\label{e3.8}
\E\sup_{t\in[0,T]}|W_{A,\sigma}(t)|^2<+\infty,
\end{equation}
we have
$$
\pi_nW_{A,\sigma}(t)\to W_{A,\sigma}(t)\quad\mbox{\rm as $n\to \infty$ in $L^2(\Omega, \mathcal F,\P),
\;\forall\;t\ge 0.$}
$$
We set $X_n(t):=X_n(t,x)$ ($=$ solution of \eqref{e3.2}). Defining
$$
W_{A_n,\sigma_n}(t)=\int_0^te^{(t-s)A_n}\sigma_n dW_n(t),\quad t\ge 0,
$$
and
$$
Z_n(t):=X_n(t)-W_{A_n,\sigma_n}(t),\quad n\in \N,\; t\ge 0,
$$
it is enough to show that
\begin{equation}
\label{e3.9}
Z_n(t)\to Z(t)\;\quad\mbox{\rm as $n\to \infty$ in $L^2(\Omega, \mathcal F,\P),\;\forall\;t\ge 0,$}
\end{equation}
because then by \eqref{e3.7}
$$
X_n(t)\to X(t)\quad\mbox{\rm as $n\to\infty$  in
$L^2(\Omega,\mathcal F,\P)$,\;$\forall \;t\ge 0$,}
$$
and  the assertion follows by Lebesgue's dominated convergence theorem. To show  \eqref{e3.9}
we first note that by the same argument as above
$$
dZ_n(t)=[A_nZ_n(t)+F_n(Z_n(t)+W_{A_n,\sigma_n}(t))]dt
$$
and thus (in the variational sense), since $A=A_n$ on $H_n$ by
\eqref{e3.1}
$$
d(Z(t)-Z_n(t))=[A(Z(t)-Z_n(t))+F_0(X(t))-F_n(X_n(t))]dt.
$$
Applying It\^o's formula we obtain that for some constant $c>0$
\begin{multline*}
\frac12\;|Z(t)-Z_n(t)|^2\le \int_0^t\big[ (\omega +1/2)|Z(s)-Z_n(s)|^2\\+|F_0(X(s))-F_0(X_n(s))|^2
+|(1-\pi_n)F_0(X(s))|^2\big]ds\\
\le c\int_0^t|Z(s)-Z_n(s)|^2ds+c\int_0^t|W_{A,\sigma}(s)-W_{A_n,\sigma_n}(s)|^2ds\\
+\int_0^t|(1-\pi_n)F_0(X(s))|^2ds.
\end{multline*}
Now \eqref{e3.9} follows by the linear growth of $F_0$,
\eqref{e3.6}-\eqref{e3.8} and Gronwall's lemma, if we can show that
\begin{equation}
\label{e3.9'}
\int_0^T\E|W_{A,\sigma}(s)-W_{A_n,\sigma_n}(s)|^2ds\to 0\quad\mbox{\rm as}\;n\to \infty.
\end{equation}
To this end we first note that a straightforward application of Duhamel's formula yields that
$$
e^{tA}|_{H_n}=e^{tA_n}\quad\;\forall\;t\ge 0.
$$
Therefore
$$
W_{A,\sigma}(s)-W_{A_n,\sigma_n}(s)=\int_0^s e^{(t-r)A}(\sigma-\pi_n\sigma\pi_n)dW(r),
$$
and thus
\begin{multline*}
\E|W_{A,\sigma}(s)-W_{A_n,\sigma_n}(s)|^2=\int_0^s \|e^{(t-r)A}(\sigma-\pi_n\sigma\pi_n)\|_{HS}^2dr\\
=\sum_{i=1}^\infty\int_0^s |e^{rA}(\sigma-\pi_n\sigma\pi_n)e_i|^2dr.
\end{multline*}
Since for any $i\in \N$, $r\in[0,s],$ the integrands converge to $0$, Lebesgue's dominated convergence
 theorem implies \eqref{e3.9'}.\bigskip

(ii) Fix $T>0$, $n\in \N$ and $x,y\in H_n$. Let $\xi^T\in C^1([0,\infty))$ be defined by
 $$
\xi^T(t):=\frac{2\omega e^{-\omega t}|x-y|}{1-e^{-2\omega T}},\quad
t\ge 0.
$$
Consider for $X_n(t)=X_n(t,x),\;t\ge 0,$ see the proof of (i), the stochastic equation
\begin{equation}
\label{e3.10} \left\{\begin{array}{l} \ds
dY_n(t)=\left[A_nY_n(t)+F_n(Y_n(t))+\xi^T(t)\;\frac{X_n(t)-Y_n(t)}{|X_n(t)-Y_n(t)|}
\;\one_{X_n(t)\neq Y_n(t)}\right]dt\\\hspace{50mm}+\sigma_ndW_n(t),\\\\
Y_n(0)=y.
\end{array}\right.
\end{equation}
Since
$$
z\to \frac{X_n(t)-z}{|X_n(t)-z|}\;\one_{X_n(t)\neq z}
$$
is dissipative on $H_n$ for all $t\ge 0$ (cf \cite{12}),
\eqref{e3.10} has a unique strong solution $Y_n(t)=Y_n(t,y),\;t\ge
0$, which is pathwise continuous $\P$-a.s.

Define the first coupling time
\begin{equation}
\label{e3.11}
\tau_n:=\inf\{t\ge 0:\;X_n(t)=Y_n(t)\}.
\end{equation}
Writing the equation for $X_n(t)-Y_n(t),\;t\ge 0,$ applying the chain rule  to
$\phi_\epsilon(z):=\sqrt{z+\epsilon^2},\;z\in (-\epsilon^2,\infty),\;\epsilon>0,$ and letting $\epsilon\to 0$ subsequently, we obtain

$$ \frac{d}{dt}\;|X_n(t)-Y_n(t)|\le
\omega\;|X_n(t)-Y_n(t)|-\\\xi^T(t)  \one_{X_n(t)\neq Y_n(t)}\quad
t\ge 0,
$$
which yields
\begin{equation}
\label{e3.12} d(e^ {-\omega t}|X_n(t)-Y_n(t)|)\le -e^{-\omega
t}\xi^T(t) \one_{X_n(t)\neq Y_n(t)}dt,\quad t\ge 0.
\end{equation}
In particular, $t\mapsto  e^{-\omega t}|X_n(t)-Y_n(t)|$ is
decreasing, hence $X_n(T)=Y_n(T)$ for all $T\ge \tau_n$. But by
\eqref{e3.12} if $T\le \tau_n$ then

$$
|X_n(T)-Y_n(T)|e^{-\omega T}\le |x-y|-|x-y|\int_0^T\frac{2\omega
e^{-2\omega t}}{1- e^{-2\omega T}}\, dt=0.$$  So, in any case
\begin{equation}
\label{e3.13}
\mbox{\rm $X_n(T)=Y_n(T)$,\quad $\P$-a.s}.
\end{equation}
Let

\begin{multline*}
R:=\exp\Bigg[-\int_0^{T\wedge
\tau_n}\frac{\xi^T(t)}{|X_n(t)-Y_n(t)|}\;\langle
X_n(t)-Y_n(t),\sigma^{-1}dW_n(t) \rangle
\\-\frac12\;\int_0^{T\wedge \tau_n}
\frac{(\xi^T(t))^2|\sigma^{-1}(X_n(t)-Y_n(t))|^2}{|X_n(t)-Y_n(t)|^2
} \;dt\Bigg]
\end{multline*}
 By \eqref{e3.13} and   Girsanov's theorem for  $p>1$,

 \begin{multline}
\label{e3.14}
(P_T^nf(y))^p= (\E[f(Y_n(T))])^p=(\E[Rf(X_n(T))])^p\\\le(P^n_Tf^p(x))(\E [R^{p/(p-1)}])^{p-1}.
\end{multline}
Let

\begin{multline*}
M_p=\exp\Bigg[-\frac{p}{p-1}\;\int_0^{T\wedge
\tau_n}\frac{\xi^T(t)}{|X_n(t)-Y_n(t)|} \;\langle
X_n(t)-Y_n(t),\sigma^{-1}dW_n(t)   \rangle
\\-\frac{p^2}{2(p-1)^2}\;\int_0^{T\wedge \tau_n}
\frac{(\xi^T(t))^2|\sigma^{-1}(X_n(t)-Y_n(t))|^2}{|X_n(t)-Y_n(t)|^2
} \;dt   \Bigg]
\end{multline*}
We have $\E M_p=1$ and hence,

\begin{multline*}
\E R^{p/(p-1)}=\E\left\{M_p\exp\left[
\frac{p}{2(p-1)^2}\;\int_0^{T\wedge \tau_n}\frac{(\xi^T(t))^2|
\sigma^{-1}(X_n(t)-Y_n(t))|^2}{|X_n(t)-Y_n(t)|^2 } \;dt  \right]
\right\}\\ \le \sup_{\Omega}\exp \left[
\frac{p}{2(p-1)^2}\;\int_0^{T\wedge \tau_n}(\xi^T(t))^2
\|\sigma^{-1}\|^2\, dt \right]\\
\le \exp\left[ \|\sigma^{-1}\|^2\frac{p\omega |x-y|^2}{(p-1)^2
(1-e^{-2\omega T})}  \right].
\end{multline*}
Combining this with \eqref{e3.14} we get the assertion (with $T$ replacing $t$). $\Box$

\section{Proof and consequences of Theorem \ref{t1.6}}
On the basis of Propositions \ref{p3.1} and \ref{p2.1} we can now easily prove Theorem \ref{t1.6}.

\noindent{\bf Proof of Theorem \ref{t1.6}}. Let $f\in Lip_b(H)$,
$f\ge 0.$ By Proposition \ref{p3.1}(i) it then follows that
\eqref{e3.4} holds with $P_tf$ replacing $P^n_tf$ provided $F$ is
Lipschitz. Using that $\bigcup_{n\in \N}H_n$ is dense in $H$ and
that $P_tf(x)$ is continuous on $x$ (cf. \cite{5}) we obtain
\eqref{e3.4} for all $x,y\in H$. In particular, this is true for
$P^{\alpha_n,\beta_n}_tf$ from Proposition \ref{p2.1}.

Now fix $t>0$ and $k\in \N$, let
$$
\chi_k(s):=\frac1k\;\one_{[t,t+1/k]}(s),\quad s\ge 0.
$$
Using \eqref{e3.4} for $P^{\alpha_n,\beta_m}_tf$, \eqref{e1.5},
Proposition \ref{p2.1} and Jensen's inequality,  we obtain for
$x,y\in K$
$$
\begin{array}{l}
\ds p_t^\mu f(x)=\lim_{k\to \infty}\frac1k\;\int_t^{t+1/k}p_s^\mu f(x)ds\\
\\
\ds=\lim_{k\to \infty}\lim_{n\to
\infty}\lim_{m\to\infty}\int_0^{t+1}
\chi_k(s)P^{\alpha_n,\beta_m}_sf(x)dx\\
\\
\ds\le\lim_{k\to \infty}\lim_{n\to
\infty}\lim_{m\to\infty}\int_0^{t+1}\chi_k(s)(P^{\alpha_n,\beta_m}_sf^p(y))^{1/p}
\exp\left[\|\sigma^{-1}\|^2\;\frac{\omega|x-y|^2}{(p-1) (1-e^{-2\omega s})}   \right]ds\\
\\
\ds\le\lim_{k\to \infty}\lim_{n\to
\infty}\lim_{m\to\infty}\left(\int_0^{t+1}\chi_k(s)P^{\alpha_n,\beta_m}_sf^p(y)
\exp\left[\|\sigma^{-1}\|^2\;\frac{p\omega|x-y|^2}{(p-1) (1-e^{-2\omega s})}  \right]ds\right)^{1/p}\\
\\
\ds=(p_t^\mu f^p(y))^{1/p}\exp\left[\|\sigma^{-1}\|^2\;\frac{\omega|x-y|^2}{(p-1)
(1-e^{-2\omega t})}   \right],
\end{array}
$$
where we note that we have to choose the sequences
$(\alpha_n),\;(\beta_n)$ such that \eqref{e2.7} holds both for $f$ and $f^p$ instead of $f$.
Since $K$ is dense in $H_0$, \eqref{e1.9} follows for $f\in C_b(H),$ for all $x,y\in H_0$, since
$p_t^\mu f$ is continuous on $H_0$ by \eqref{e1.3}.

Let now $f\in B_b(H),\;f\ge 0$. Let $f_n\in C_b(H)$, $n\in \N$, such that $f_n\to f$ in $L^p(H,\mu)$ as
$n\to \infty,$ $p\in (1,\infty)$ fixed. Then, since $\mu$ is invariant for $p_t^\mu,\;t>0, $ selecting a
subsequence if necessary, it follows that there exists  $K_1\in \mathcal B(H),$ $\mu(K_1)=1$, such that
$$
p_t^\mu f_n(x)\to p_t^\mu f(x)\quad\mbox{\rm as}\;n\to \infty,\;\forall\; x\in K_1.
$$
Taking this limit in \eqref{e1.9} we obtain \eqref{e1.9} for all $ x,y\in K_1.$ Taking into account that
$p_t^\mu$ is continuous and that $K_1$ is dense in $H_0=$ supp $\mu$, \eqref{e1.9}  follows for all
$x,y\in H_0$. $\Box$ \bigskip

Corollary \ref{c1.7} immediately follows from Theorem \ref{t1.6} and the following general result:
\begin{Proposition}
\label{p4.1}
 Let $E$ be a topological space and $P$ a Markov operator on $B_b(E)$. Assume that for any $p>1$ there exists
 a continuous function $\eta_p$ on $E\times E$ such that $\eta_p(x,x)=0$ for all $x\in E$ and
\begin{equation}
\label{e4.1}
P|f|(x)\le (P|f|^p(y))^{1/p} e^{\eta_p(x,y)}\quad\;\forall\;x,y\in E,\;f\in B_b(E).
\end{equation}
Then $P$ is strong Feller, i.e. maps $B_b(E)$ into $C_b(E)$. Furthermore, for any $\sigma$-finite measure $\mu$ on
$(E,\mathcal B(E))$ such that
\begin{equation}
\label{e4.2}
\int_E |Pf|d\mu \le C\int_E |f|d\mu,\quad\;\forall\;
f\in B_b(E),
\end{equation}
for some $C>0$, $P$ uniquely extends to $L^p(E,\mu)$ with
$PL^p(E,\mu)\subset C(E)$ for any $p> 1$.
\end{Proposition}
{\bf Proof}. Since $P$ is linear, we only need to consider $f\ge 0$. Let $f\in B_b(E)$ be nonnegative.
By \eqref{e4.1} and the property of $\eta_p$ we have
$$
\limsup_{x\to y}Pf(x)\le(Pf^p(y))^{1/p},\quad p>1.
$$
Letting $p\downarrow 1$ we obtain $\limsup_{x\to y}Pf(x)\le Pf(y)$. Similarly, using $f^{1/p}$ to replace $f$
and replacing $x$ with $y$, we obtain
$$
(Pf^{1/p}(y))^p\le (Pf(x))e^{p\eta_p(y,x)},\quad\;\forall\;
x,y\in E,\; p>1.
$$
First letting $x\to y$ then $p\to 1$, we obtain $\liminf_{x\to
y}Pf(x)\ge Pf(y)$. So $Pf\in C_b(E)$. Next, for any nonnegative
$f\in L^p(E,\mu)$, let $f_n=f\wedge n,\; n\ge 1$. By \eqref{e4.2}
  and $f_n\to f$ in $L^p(E,\mu)$ we
have $P|f_n-f_m|^p\to 0$ in $L^1(E,\mu)$ as $n, m\to \infty$. In
particular, there exists $y\in E$ such that
\begin{equation}
\label{e4.3}
\lim_{n, m\to \infty}P|f_n-f_m|^p(y)=0.
\end{equation}
Moreover, by \eqref{e4.1}, for $B_N:=\{x\in E:\; \eta_p(x,y)< N\}$
$$
\sup_{x\in B_N}|Pf_n(x)-Pf_m(x)|^p \le \sup_{x\in B_N}(P|f_n-f_m|(x))^p\le (P|f_n-f_m|^p(y))e^{pN}.
$$
Since by the strong Feller property $Pf_n\in C_b(E)$ for any $n\ge 1$ and noting that $C_b(B_N)$ is complete
under the  uniform norm, we conclude from \eqref{e4.3}
that $Pf$ is continuous on $B_N$ for any $N\ge 1$, and hence, $Pf\in C(H).$ $\Box$\bigskip

\

\noindent{\bf Proof of Corollary \ref{c1.8}}. Let $\mu_1,\mu_2$ be
probability measures on $(H,\mathcal B(H))$ satisfying (H4). Define
$\mu:=\frac12\;\mu_1+\frac12\;\mu_2$. Then $\mu$ satisfies (H4) and
$\mu_i=\rho_i\mu,\;i=1,2,$ for some $\mathcal B(H)$-measurable
$\rho_i:H\to [0,2].$ Let $i\in \{1,2\}.$

Since $\rho_i$ is bounded, by(H4)(iii) and Theorem \ref{t1.2} it follows that
$$
\int_H L_\mu u\;d\mu_i=0,\quad\;\forall\; u\in D(L_\mu).
$$
Hence
$$
\frac{d}{dt}\;\int_H e^{tL_\mu} u\;d\mu_i=
\int_HL_\mu(e^{tL_\mu} u)d\mu_i=0,\quad\;\forall\; u\in D(L_\mu),
$$
i.e.
$$
\int_H p_t^\mu u\;d\mu_i=\int_H u\;d\mu_i\quad\;\forall\; u\in \mathcal E_A(H).
$$
Since $\mathcal E_A(H)$ is dense in $L^1(H,\mu_i)$, $\mu_i$ is $(p_t^\mu)$-invariant. But as mentioned before,
by Theorem \ref{t1.6} it follows that $(p_t^\mu)$ is irreducible on $H_0$ (see \cite{9a}) and it is strong Feller
on $H_0$ by Corollary \ref{c1.7}. So, since $\mu_i(H_0)=1$, $\mu_i=\mu$. $\Box$\bigskip

\noindent{\bf Proof of Corollary \ref{c1.10}}.
Let
$$
\begin{array}{l}
\tilde A:=A-\omega I,\quad D(\tilde A):=D(A)
\\
\tilde F_0:=F_0+\omega I.
\end{array}
$$
By $(H2)$, $\tilde A$ has discrete spectrum. Let $e_k\in H$, $-\lambda_k\in (-\infty,0]$, be the corresponding orthonormal eigenvectors, eigenvalues respectively.

For $k\in \N$ define
$$
\varphi_k(x):=\langle e_k,x  \rangle,\quad x\in H.
$$
We note that by a simple approximation \eqref{e1.4} also holds for any Lipschitz function on $H$ and thus (cf. the proof of \cite[Proposition 5.7(iii)]{4})
also \eqref{e1.5} holds for such functions, i.e. in particular, for all $k\in \N$
\begin{equation}
\label{e4.4ro}
[0,\infty)\ni t\mapsto p_t\varphi_k(x)\quad\mbox{\rm is continuous for all $x\in H_0$.}
\end{equation}
Since any compactly supported smooth function on $\mathbb R^N$ is
the Fourier
       transform of a Schwartz test function, by approximation it easily   follows that setting
$$
\mathcal F C^\infty_b(\{e_k\}):=\{g(\langle e_1,\cdot \rangle,...,\langle e_N,\cdot \rangle) :\;N\in \N,\;g\in C^\infty_b(\R^N)\},
$$
we have $\mathcal F C^\infty_b(\{e_k\})\subset D(L_\mu)$ and for $\varphi\in \mathcal F  C^\infty_b(\{e_k\})$
$$
L_\mu\varphi(x)=\frac12\;\mbox{\rm Tr}\;[D^2\varphi(x)]+\langle x, AD\varphi(x) \rangle+\langle F_0(x), D\varphi(x) \rangle
\quad x\in H.
$$
Then by approximation
it is easy to show that
\begin{multline}
\label{e4.5ro}
\varphi_k, \varphi^2_k\in D(L_\mu)\;\mbox{\rm and}\; L_\mu\varphi_k=-\lambda_k\varphi_k+\langle e_k,\tilde F_0 \rangle,
\\ L_\mu\varphi^2_k=-2\lambda_k\varphi^2_k+2\varphi_k\langle e_k,\tilde F_0 \rangle+2
 \quad\forall\;k\in \N.
\end{multline}
Since we assume that $|F_0|$ is in $L^2(H,\mu)$, by \cite[Theorem
1.1]{4'}
       we are in the situation of \cite[Chapter II]{13'}. So, we conclude that by
 \cite[Chapter II, Theorem 1.9]{13'} there exists a normal (that is $\P_x[X(0)=x]=1$) Markov process
$(\Omega,\mathcal F,(\mathcal F_t)_{t\ge 0},(X(t))_{t\ge
0},(\P_x)_{x\in H_0})$ with state space $H_0$ and $M\in \mathcal
B(H_0)$,   $\mu(M)=1$, such that $X(t)\in M$ for all $t\ge 0$
$\P_x$-a.s. for all $x\in M$ and which has continuous sample paths
$\P_x$-a.s for all $x\in M$ and for which by the proof of
\cite[Proposition 8.2]{4} and \eqref{e4.4ro}, \eqref{e4.5ro} we have
that for all $k\in \N$
\begin{equation}
\label{e4.5ro'}
\begin{array}{l}
\ds\beta^x_k(t):=\varphi_k(X(t))-\varphi_k(x)-\int_0^tL_\mu\varphi_k(X(s))ds,\quad t\ge 0,\\
\\
\ds
M^x_k(t):=\varphi^2_k(X(t))-\varphi^2_k(x)-\int_0^tL_\mu\varphi^2_k(X(s))ds,\quad
t\ge 0,
\end{array}
\end{equation}
are continuous   local $(\mathcal F_t)$-martingales with $\beta^x_k(0)=M_k(0)=0$ under $\P_x$ for all $x\in M$. Fix $x\in M$.
Below $\E_x$ denotes expectation with respect to $\P_x$. Since for $T>0$
$$
\begin{array}{l}
\ds\int_H\int_0^T\E_x(1+|X(s)|^2)(1+|F_0(X(s))|) ds\mu(dx)\\
\\
\ds=T
\int_H(1+|x|^2)(1+|F_0(x)|)\mu(dx)<\infty,
$$
\end{array}
$$
making $M$ smaller if necessary,
by (H4)(ii) we may assume that
\begin{equation}
\label{e4.6ro}
\E_x\int_0^T(1+|X(s)|^2)(1+|F_0(X(s))|)ds<\infty.
\end{equation}
 By standard Markov process theory we have for their    covariation   processes under $\P_x$,
\begin{equation}
\label{e4.6ro'}
\langle\beta^x_k,\beta^x_{k'}   \rangle_t=\int_0^t\langle D\varphi_k(X(s)) , D\varphi_{k'}(X(s))\rangle  ds=t\delta_{k,k'},\quad t\ge 0.
\end{equation}
Indeed, an
elementary calculation shows that for all $k\in\N,\;t\ge 0,$
\begin{multline}
\label{e4.6ro''}
\beta^x_k(t)^2-\int_0^t|D\varphi_k(X(s))|^2ds\\=M_k^x(t)-2\varphi_k(x)\beta^x_k(t)-\int_0^t(\beta^x_k(t)-\beta^x_k(s))L_\mu\varphi_k(X(s)) ds,
\end{multline}
where all three summands on the right hand side are martingales. Since we have a similar formula for finite linear combinations of $\varphi_k's$ replacing a single $\varphi_k$, by polarization we get \eqref{e4.6ro'}. Note that by \eqref{e4.5ro} and \eqref{e4.6ro} all integrals in \eqref{e4.5ro'}, \eqref{e4.6ro''} are well defined.

Hence, by \eqref{e4.6ro'} $\beta^x_k,\;k\in \N,$ are independent standard $(\mathcal F_t)$-Brownian motions under $\P_x$. Now it follows by \cite[Theorem 13]{11'} that, with
$W^x=(W^x(t))_{t\ge 0},$ being the cylindrical Wiener process on $H$
given by $W^x=(\beta^x_ke_k)_{k\in \N},$ we have for every $t\ge 0$,
\begin{equation}
\label{e4.7ro}
X(t)=e^{tA}x+\int_0^te^{(t-s)A}F_0(X(s))ds+\int_0^te^{(t-s)A}
dW^x(s),\quad\P\mbox{\rm -a.s.},
\end{equation}
that is, the tuple $(\Omega,\mathcal F,(\mathcal F_t)_{t\ge 0},\P_x,W^x,X)$
is a solution to
\begin{equation}
\label{e4.8ro}
\left\{\begin{array}{l}
\ds Y(t)=e^{tA}Y(0)+\int_0^te^{(t-s)A}F_0(Y(s))ds+\int_0^te^{(t-s)A}  dW(s),\quad\P\mbox{\rm -a.s.},\quad\forall\;t\ge 0,\\
\mbox{\rm law  $Y(0)=\delta_x (:=$ Dirac measure in $x$)},
\end{array}\right.
\end{equation}
in the sense of \cite[page 4]{11'}.

We note that the zero set in \eqref{e4.7ro} is indeed independent of $t$, since all terms are continuous in $t$ $\P_x$-a.s. because of (H2)(ii)
and \eqref{e4.6ro}.\bigskip

\noindent{\bf Claim} We have $X$-pathwise uniqueness
for equation \eqref{e4.8ro} (in the sense of \cite[page 98]{11'}).\bigskip

For any given cylindrical $(\mathcal F'_t)$-Wiener process $W$ on a stochastic basis $(\Omega',\mathcal F',(\mathcal F'_t)_{t\ge 0},\P')$ let $Y=Y(t),\;Z=Z(t),\;t\ge 0,$ be two solutions of \eqref{e4.8ro} such that law$(Z)$=law$(Y)$=law$(X)$ and $Y(0)=Z(0)\;\P'$-a.s.. Then by \eqref{e4.6ro}
\begin{equation}
\label{e4.9ro}
\E '\int_0^T|F_0(Y(s))|ds=\E'\int_0^T|F_0(Z(s))|ds=\E_x\int_0^T|F_0(X(s))|ds<\infty.
\end{equation}
(which, in particular implies by \eqref{e4.8ro} and by (H2)(i) that
both $Y$ and $Z$ have $\P'$-a.s. continuous sample paths). Hence
applying \cite[Theorem 13]{11'} again (but this time using the dual
implication) we obtain for all $k\in \N$
$$
\begin{array}{lll}
 \langle e_k, Y(t)-Z(t)  \rangle&=&\ds-\lambda_k\int_0^t\langle  e_k, Y(s)-Z(s)  \rangle ds\\
 \\
 &&\ds+\int_0^t\langle  e_k, \tilde{F_0}(Y(s))-\tilde{F_0}(Z(s)) \rangle ds,\quad t\ge 0,\;\P'\mbox{\rm -a.s.}.
 \end{array}
$$
Therefore, by the chain rule for all $k\in \N$
\begin{multline*}
\langle e_k, Y(t)-Z(t)  \rangle^2=-2\lambda_k\int_0^t\langle  e_k, Y(s)-Z(s)  \rangle^2ds\\
\\
+2\int_0^t\langle  e_k, Y(s)-Z(s)  \rangle \;\langle  e_k, \tilde{F_0}(Y(s))-\tilde{F_0}(Z(s)) \rangle ds,\quad t\ge 0,\;\P'\mbox{\rm -a.s.}.
\end{multline*}
Dropping the first term on the right hand side  and summing up over
$k\in \N$ (which is justified by \eqref{e4.8ro} and the continuity
of $Y$ and $Z$), we obtain from (H3) that
\begin{multline*}
|Y(t)-Z(t)|^2\le 2\int_0^t\langle   Y(s)-Z(s), \tilde{F_0}(Y(s))-\tilde{F_0}(Z(s)) \rangle ds\\
\\
\le 2\omega\int_0^t    |Y(s)-Z(s)|^2ds,\quad t\ge 0,\;\P'\mbox{\rm -a.s.}.
\end{multline*}
Hence, by Gronwall's lemma $Y=Z$ $\P'$-a.s. and the Claim is proved.

 By the Claim we can apply \cite[Theorem 10, $(1)\Leftrightarrow (3)$]{11'} and then \cite[Theorem 1]{11'} to conclude that equation \eqref{e4.8ro} has a strong solution (see \cite[Definition 1]{11'}) and that there is one strong solution with the same law as $X$, which hence by \eqref{e4.6ro} has continuous sample paths a.s. Now all conditions in \cite[Theorem 13.2]{11'} are fulfilled and, therefore, we deduce from it that on any stochastic basis $(\Omega,\mathcal F,(\mathcal F_t)_{t\ge 0},\P)$ with $(\mathcal F_t)$-cylindrical Wiener process $W$ on $H$ and for $x,y\in M$ there exist pathwise unique continuous strong solutions $X(t,x),\;X(t,y),\;t\ge 0,$ to \eqref{e4.8ro} such that
 $$
 \P\circ X(\cdot,x)^{-1}=\P_x\circ X^{-1}
 $$
and
$$
 \P\circ X(\cdot,y)^{-1}=\P_y\circ X^{-1},
 $$
in particular, $X(0,x)=x$ and $X(0,y)=y$ and
\begin{equation}
\label{e4.10ro}
\begin{array}{l}
\P\circ X(t,x)^{-1}(dz)=p_t(x,dz),\quad t\ge 0,\\
\P\circ X(t,y)^{-1}(dz)=p_t(y,dz),\quad t\ge 0.
\end{array}
\end{equation}
In particular, we have proved (i). To prove (ii),
below for brevity we set $X:=X(\cdot,x),\;X':=X(\cdot,y)$. Then proceeding as in the proof of the Claim, by \eqref{L2} and noting
that $s^{-1}\Phi(s)\to\infty$ as
$s\to\infty,$ we obtain
\begin{equation}
\label{e4.5r}
 \frac{d}{dt}\; |X(t)-X'(t)|^2\le a - \Phi_0(|X(t)-X'(t)|^2)
\end{equation}
  for some constant $a>0$, only depending on $\omega$ and $\Phi$, where $\Phi_0=\frac12\;\Phi$.

 Now we  consider two cases.\bigskip

{\bf Case 1}. $|x-y|^2\le \Phi_0^{-1}(2a)$.\bigskip

Define $f(t):=|X(t)-X'(t)|^2,\;t\ge 0,$ and suppose there exists $t_0\in (0,\infty)$ such that
$$
f(t_0)>\Phi_0^{-1}(a).
$$
Then we can choose $\delta\in [0,t_0]$ maximal such that
$$
f(t)>\Phi_0^{-1}(a),\quad\forall\;t\in(t_0-\delta,t_0].
$$
Hence, because  by \eqref{e4.5r} $f$ is decreasing on every interval where it is larger than $\Phi_0^{-1}(a)$, we obtain that
$$
f(t_0-\delta)\ge f(t_0)>\Phi_0^{-1}(a).
$$
Suppose $t_0-\delta>0$. Then $f(t_0-\delta)\le \Phi_0^{-1}(a)$ by the continuity of $f$ and  the maximality of $\delta$. So, we  must have $t_0-\delta=0$, hence
$$
f(t_0)\le f(t_0-\delta)=f(0)=|x-y|^2\le \Phi_0^{-1}(2a).
$$
So,
$$|X(t)-X'(t)|^2\le \Phi_0^{-1}(2a),\quad\forall\;t>0.$$

{\bf Case 2}. $|x-y|^2> \Phi_0^{-1}(2a)$.\bigskip

 Define $t_0=\inf\{t\ge 0: |X(t)-X'(t)|^2\le \Phi_0^{-1}(2a)\}.$  Then by Case 1, starting at $t=t_0$ rather than $t=0$ we know that
 \begin{equation}
\label{e4.6r}
  |X(t)-X'(t)|^2\le  \Phi_0^{-1}(2a),\quad\forall\;t\ge t_0.
\end{equation}
 Furthermore, it follows
from \eqref{e4.5r} that
$$d |X(t)-X'(t)|^2\le -\ff 1 2\Phi_0(|X(t)-X'(t)|^2)\d t,\quad \forall\; t\le
t_0.$$
 This implies
$$\Psi(|X(t)-X'(t)|^2)\ge \frac12\;\int_{|X(t)-X'(t)|^2}^{|x-y|^2}\ff{d
r}{\Phi_0(r)}\ge \ff t 4,\quad \forall\; t\le t_0.$$
 Therefore,
 \begin{equation}
\label{e4.7r}
 |X(t)-X'(t)|^2\le \Psi^{-1}(t/4),\quad \forall\; t\le t_0.
\end{equation}
Combining Case 1, \eqref{e4.6r} and  \eqref{e4.7r}  we conclude that
\begin{equation}
\label{L4}
|X(t)-X'(t)|^2\le
\Psi^{-1}(t/4)+ \Phi_0^{-1}(2a),\quad \forall\;t>0.
\end{equation}
Combining \eqref{L4} with
Theorem \ref{t1.6} for all $f\in B_b(H)$ we obtain
$$\big(p_{t/2} |f|(X(t/2))\big)^2\le
\big(p_{t/2}f^2(X'(t/2)\big)\exp\Big[\ff{\ll(1+\Psi^{-1}(t/8))}{(1-\e^{-\oo
t/2})^2}\Big],\quad \forall\;t>0$$ for some constant $\ll>0.$
 By Jensen's inequality and approximation it follows that for all $f\in L^2(H,\mu)$
\begin{multline}
\label{4.14ro}
(p_t |f|(x))^2\le \mathbb E\big(p_{t/2} |f|(X(t/2))\big)^2\\\le
\big(p_t
f^2(y)\big)\exp\Big[\ff{\ll(1+\Psi^{-1}(t/8))}{(1-\e^{-\oo
t/2})^2}\Big],\quad \forall\; t>0,\;\forall\;x,y\in M.
\end{multline}
But since $H_0=$ supp $\mu$, $M$ is dense in $H_0$, hence by the continuity of $p_tf$ (cf. Corollary \ref{c1.7}) \eqref{4.14ro} holds for all $x\in H_0$, $y\in M$. Since $\mu(M)=1$ this
completes the proof by integrating
both sides with respect to $\mu(dy).$
$\Box$
\begin{Remark}
\label{r4.2ro} \em We would like to  mention that by using
\cite{BBR} instead of \cite{13'} we can drop
       the assumption that $|F_0| \in L^2(H,\mu).$ So, by (\ref{e4.6ro''}) and the proof
       above we can derive (\ref{e4.6ro'}) avoiding
 to assume the usually energy condition
$$
\int_0^t|F_0(X(s))|^2ds<\infty,\quad \P_x\mbox{\rm -a.s.}.
$$
Details will be included in a forthcoming
        paper.
 We would like to thank Tobias Kuna at this point from whom we learnt identity \eqref{e4.6ro''} by private communication.
\end{Remark}

\section{Existence of measures satisfying (H4)}

To prove existence of invariant measures we need to strengthen some of our assumptions. So, let us introduce the
following conditions.\bigskip

\noindent (H1)' $(A,D(A))$ is self-adjoint satisfying \eqref{e1.1'}.

\noindent (H6) There exists $\eta\in (\omega,\infty)$ such that
$$
\langle F_0(x)-F_0(y),x-y   \rangle
\le -\eta|x-y|^2,\quad\;\forall\;
x,y\in D(F).
$$
\begin{Remark}
\label{r5.1}
\em (i) Clearly, (H1)' implies (H1) and (H5). (H1)' and  (H2)(i) imply that $(A,D(A))$ and thus also
$(1+\omega-A,D(A))$ has a discrete spectrum. Let $\lambda_i\in (0,\infty),\;i\in \N,$ be the eigenvalues of the
latter operator. Then by (H2)
\begin{equation}
\label{e5.1}
\sum_{i=1}^\infty\lambda_i^{-1}<\infty.
\end{equation}

\noindent (ii) If we assume \eqref{e5.1}, i.e. that  $(1+\omega-A)^{-1}$ is trace class, then all what follows
holds with (H2) replaced by (H2)(i). So, $\sigma^{-1}\in L(H)$ is not needed in this case.
\end{Remark}

Let $F_\alpha$, $\alpha<0,$ be as in Section 2. Then e.g. by
\cite[Theorem 3.2]{DP} equation \eqref{e1.1} with $F_\alpha$
replacing $F_0$ has a unique mild solution $X_\alpha(t,x),\;t\ge 0.$
Since there exist $\tilde \eta\in (\omega,\infty)$ and $\alpha_0>0$
such that each $F_\alpha$, $\alpha\in (0,\alpha_0)$, satisfies (H6)
with $\tilde \eta$ replacing $\eta$, by \cite[Section 3.4]{DP}
$X_\alpha$ has a unique invariant measure $\mu_\alpha$ on
$(H,\mathcal B(H))$ such that for each $m\in \N$
\begin{equation}
\label{e5.2}
\sup_{\alpha\in (0,\alpha_0)}\int_H|x|^m\mu_\alpha(dx)<\infty.
\end{equation}
That these moments are indeed uniformly bounded  in $\alpha$,
follows from the proof of \cite[Proposition 3.18]{DP} and the  fact
that $\tilde \eta\in (\omega,\infty)$.

Let $N_Q$ denote the centered  Gaussian measure on $(H,\mathcal B(H))$ with covariance operator $Q$ defined by
$$
Qx:=\int_0^\infty e^{tA}\sigma e^{tA}xdt,\quad x\in H,
$$
which by (H2)(ii) is  trace class.

Let $W^{1,2}(H,N_Q)$ be defined as usual, that is as the completion of $\mathcal E_A(H)$ with respect to the norm
$$
\|\varphi\|_{W^{1,2}}:=\left(\int_H(\varphi^2+|D\varphi|^2) dN_Q  \right)^{1/2},\quad \varphi\in \mathcal E_A(H),
$$
where $D$ denotes first Fr\'echet derivative. By \cite{4''} we know that
\begin{equation}
\label{e5.3}
W^{1,2}(H,N_Q)\subset L^2(H,N_Q),\quad\mbox{\rm compactly.}
\end{equation}

\begin{Theorem}
\label{t5.2}
Assume that (H1)', (H2), (H3) and (H6) hold and let $\mu_\alpha, \alpha\in (0,\alpha_0)$ be as above. Suppose that there exists a lower semi-continuous
function $G:H\to [0,\infty]$ such that
\begin{equation}
\label{e5.4}
\{G<\infty\}\subset D(F),\quad |F_0|\le G\;\mbox{\rm on $D(F)$ and}\;
\sup_{\alpha\in (0,\alpha_0)}\int_H G^2d\mu_\alpha<\infty.
\end{equation}
Then $\{\mu_\alpha:\;\alpha\in (0,\alpha_0)\}$ is tight and any
limit point $\mu$ satisfies (H4) and hence by Corollary \ref{c1.8}
all of these limit points coincide. Furthermore, for all $m\in \N$
\begin{equation}
\label{e5.5}
\int_H(|F_0(x)|^2+|x|^m)\mu(dx)<\infty
\end{equation}
and there exists $\rho:H\to[0,\infty),$ $\mathcal B(H)$-measurable, such that $\mu=\rho N_Q$ and $\sqrt\rho\in W^{1,2}(H,\mu)$.
\end{Theorem}
{\bf Proof}. We recall that by \cite[Theorem  1.1]{4'} for each $
\alpha\in (0,\alpha_0)$
\begin{equation}
\label{e5.6}
\mu_\alpha=\rho_\alpha N_Q;\quad \sqrt{\rho_\alpha}\in W^{1,2}(H,N_Q)
\end{equation}
and as is easily seen from its proof, that
\begin{equation}
\label{e5.7}
\int_H|D\sqrt{\rho_\alpha}|^2dN_Q\le \frac14\;\int_H|F_\alpha|^2d\mu_\alpha.
\end{equation}
But by \eqref{e2.3} and \eqref{e5.4} the right hand side of \eqref{e5.7} is uniformly bounded in $\alpha$.
Hence by \eqref{e5.3} there exists a zero sequence $\{\alpha_n\}$ such that
$$
\sqrt{\rho_{\alpha_n}}\to\sqrt{\rho}\quad\mbox{\rm in $L^2(H,N_Q)$ as $n\to \infty$},
$$
for some $\sqrt{\rho}\in W^{1,2}(H,N_Q)$ and therefore, in particular,
\begin{equation}
\label{e5.8}
\rho_{\alpha_n}\to \rho\quad\mbox{\rm in\;$L^1(H,N_Q)$ as $n\to \infty$}.
\end{equation}
Define $\mu:=\rho N_Q$ and $\rho_n:=\rho_{\alpha_n},\;n\in \N.$ Since $G$ is lower semi-continuous  and $\mu_{\alpha_n}\to \mu$ as $n\to \infty$ weakly,
\eqref{e5.2} and \eqref{e5.4} imply
\begin{equation}
\label{e5.9}
\int_H(G^2(x)+|x|^m)\mu(dx)<\infty\quad\forall\;m\in \N.
\end{equation}
Hence by \eqref{e5.4} both (H4)(i) and (H4)(ii) follow. So, it remains to  prove (H4)(iii).

Since $\sigma$ is independent of $\alpha$, to show \eqref{e5.9} it is  enough  to prove that for all $\varphi\in C_b(H),\;h\in D(A)$
\begin{equation}
\label{e5.10}
\lim_{n\to \infty}\int_HF^h_{\alpha_n}(x)\varphi(x)\mu_{\alpha_n}(dx)=\int_HF^h_{0}(x)\varphi(x)\mu(dx),
\end{equation}
where $F^h_{\alpha}:=\langle h,F_\alpha \rangle$, $\alpha\in [0,\alpha_0)$. We have
\begin{multline}
\label{e5.11}
\left|\int_HF^h_{\alpha_n}\varphi d\mu_{\alpha_n}-  \int_HF^h_{0}\varphi d\mu\right|\\
\le \|\varphi\|_\infty\int_H|F^h_{\alpha_n}-F_0^h|\rho_ndN_Q+\int_H|F_0^h\varphi|\;|\rho_n-\rho| dN_Q.
\end{multline}
But by \eqref{e2.3} and  \eqref{e5.4} we have
\begin{multline*}
\int_H|F^h_{\alpha_n}-F_0^h|\rho_ndN_Q\le
\int_{\{|G|\le M\}}|F^h_{\alpha_n}-F_0^h|\rho_ndN_Q\\
+\frac{2|h|}{M}\;\sup_{\alpha\in (0,\alpha_0)}\int_H
G^2 d\mu_{\alpha}.
\end{multline*}
Hence first letting $n\to  \infty$ then $M\to  \infty$
by \eqref{e2.2}, \eqref{e5.4}  and  \eqref{e5.8} Lebesgue's generalized dominated convergence theorem implies that the first term on the right hand side of \eqref{e5.11} converges to $0$. Furthermore, for every $\delta\in (0,1)$
\begin{multline}
\label{e5.12}
\left| \int_HF^h_{0}\varphi d\mu_{\alpha_n}-  \int_HF^h_{0}\varphi d\mu\right|\\
\le \left| \int_H\frac{F^h_{0}}{1+\delta|F^h_{0} |}\;\varphi (\rho_n-\rho)dN_Q  \right|\\
+\delta\|\varphi\|_\infty\left( \int_H|F^h_{0}|^2 d\mu_{\alpha_n}+ \int_H|F^h_{0}|^2  d\mu   \right).
\end{multline}
Since by \eqref{e2.3} and \eqref{e5.4}
$$
\sup_{\alpha\in (0,\alpha_0)}\int_H
 |F^h_{0}|^2 d\mu_{\alpha}<\infty,
$$
(H4)(iii) follows  from  \eqref{e5.12} by letting first $n\to \infty$ and then $\delta\to 0$, since for fixed $\delta>0$ the  first term in the right hand side converges to zero by \eqref{e5.8}. $\Box$
\begin{Example}
\label{ex5.4}
\em Let $H=L^2(0,1), Ax=\Delta x$, $\; x\in D(A):=H^2(0,1)\cap H^1_0(0,1).$   Let $f:\R\to \R$ be decreasing such  that for some $c_3>0,\; m\in \N$,
\begin{equation}
\label{e5.13}
|f(s)|\le c_3(1+|s|^m),\quad\forall\;s\in  \R.
\end{equation}
Let $s_i\in \R$, $i\in \N$, be the set of all arguments where $f$ is not continuous and define
$$
\bar f(s)=
\left\{\begin{array}{l}
\protect
[f(s_{i^+}),f(s_{i^-})],\quad\mbox{\rm if $s=s_i$ for some $i\in \N$},\\
f(s),\quad\mbox{\rm else}
\end{array}\right.
$$
Define
$$
F:D(F)\subset H\to 2^H,\;\;x\mapsto \bar f\circ x,
$$
where
$$
D(F)=\{x\in H:\;\bar f\circ x\subset H\}.
$$
Then $F$ is $m$-dissipative. Let $F_0$ be defined as in Section 2.

Since $A\le \omega$ for some $\omega<0$, it is easy to check that all conditions (H1)',  (H2), (H3), (H6) with $\eta=0$  hold for any $\sigma\in L(H)$
such that $\sigma^{-1}\in L(H)$. Define
$$
G(x):=\left\{\begin{array}{l}
\ds\left(\int_0^1|x(\xi)|^{2m}d\xi  \right)^{1/2}\quad\mbox{\rm if}\;x\in L^{2m}(0,1)\\
\\
+\infty\quad\mbox{\rm if}\;x\notin L^{2m}(0,1).
\end{array}\right.
$$
Then $\{G<\infty\}\subset D(F)$ and $|F_0|=|F_0|_{L^2(0,1)}\le G$ on $D(F)$. Furthermore, by \cite[(9.3)]{4}
\begin{equation}
\label{e5.14}
\sup_{\alpha\in (0,\alpha_0)}\int_HG^2d\mu_\alpha<\infty.
\end{equation}
Note that from \cite[Hypothesis 9.5]{4}  only the first inequality,
which clearly holds by \eqref{e5.13} in our case, was used to prove
\cite[(9.3)]{4}. Hence all assumptions of Theorem \ref{t5.2} above
hold and we obtain the existence of the desired unique probability
measure $\mu$ satisfying (H4) in this case. We emphasize that no
continuity properties of $f$ and $F_0$ are required. In particular,
then all results stated in Section 1 except for Corollary 1.10(ii)
hold in this case.

If moreover there exists an increasing positive convex function $\Phi$ on $[0,\infty)$ satisfying (\ref{L1})
such that

$$(f(s)-f(t))(s-t)\le c-\Phi(|s-t|^2),\ \ \ s,t\in \mathbb R,$$ then by  Jensen's inequality (\ref{L2}) holds. Hence,
by Corollary \ref{c1.10} one obtains an explicit upper bound for
$\|p_t\|_{L^2(H,\mu)\to L^\infty(H,\mu)}.$ A natural and simple choice of $\Phi$ is
 $\Phi(s)=s^m$ for $m>1$.

One can extend these results to  the case, where $(0,1)$ above is replaced by a  bounded open set in $\R^d$, $d=2$ or $3$ for $\sigma=(-\Delta)^\gamma$, $\gamma\in (\frac{d-2}{4},\frac12)$,
 based on Remark \ref{r1.1}(iv).
  \end{Example}\bigskip

Before to conclude we want to present a condition in the general case (i.e for any Hilbert space $H$ as above) that implies \eqref{e5.4}, hence by Theorem \ref{t5.2} ensures the existence of a probability measure satisfying (H4) so that all results of Section 1 apply also to this case. As will become clear from the arguments below, such condition is satisfied if the eigenvalues of $A$ grow fast enough in comparison with $|F_0|$. To this end we first note that by \eqref{e5.1}
  for $i\in \N$ we can find $q_i\in (0,\lambda_i)$, $q_i\uparrow \infty$ such that
 $\sum_{i=1}^\infty q_i^{-1}<\infty$ and $\frac{q_i}{\lambda_i}\to 0$ as $i\to \infty$. Define
 $\Theta:H\to [0,\infty]$ by
\begin{equation}
\label{e5.15}
\Theta(x):=\sum_{i=1}^\infty\frac{\lambda_i}{q_i}\; \langle x,e_i   \rangle^2,\quad x\in H,
\end{equation}
where $\{e_i\}_{i\in N}$ is an eigenbasis  of $(1+\omega-A,D(A))$ such that $e_i$ has eigenvalue $\lambda_i$.
Then $\Theta$  has compact level sets
and $|\cdot|^2\le \Theta$.\bigskip

 Below we set
$$
\mbox{\rm $H_n$:=lin span $\{e_1,...,e_n\}$},\quad
\mbox{\rm $\pi_n$:=projection onto $H_n$},
$$
\begin{equation}
\label{e5.16}
\tilde A:=A-(1+\omega)I,\quad D(\tilde A):=D(A),
\end{equation}
\begin{equation}
\label{e5.17}
\tilde F_0:=F_0+(1+\omega)I.
\end{equation}
We note that obviously $H_n\subset\{\Theta<+\infty\}$ for all $n\in \N$.
\begin{Theorem}
\label{t5.4}
Assume that (H1)', (H2), (H3) and (H6) hold and let $\mu_\alpha,\;\alpha\in (0,\alpha_0)$, be as above.
Suppose that $\{\Theta<+\infty\}\subset D(F)$ and that for  some $C\in(0,\infty), m\in \N$
\begin{equation}
\label{e5.18}
|F_0(x)|\le C(1+|x|^m+\Theta^{1/2}(x)),\quad\forall\;x\in D(F).
\end{equation}
Then
\begin{equation}
\label{e5.19}
\sup_{\alpha\in (0,\alpha_0)}\int_H\Theta d\mu_\alpha<\infty
\end{equation}
and \eqref{e5.4} holds, so Theorem \ref{t5.2} applies.
\end{Theorem}
{\bf Proof}. Consider the Kolmogorov operator $L_\alpha$ corresponding to $X_\alpha(t,x),\;t\ge 0,\;x\in H,$ which for $\varphi\in \mathcal FC^2_b(\{e_n\})$, i.e.,  $\varphi=g(\langle e_1,\cdot   \rangle
,...,\langle e_N,\cdot   \rangle
)$ for some  $N\in \N$, $g\in C^2_b(\R^N)$, is given by
\begin{equation}
\label{e5.20}
L_\alpha\varphi(x):=\frac12\;\mbox{\rm Tr}\;[\sigma^2 D^2\varphi(x)]+\langle x,AD\varphi(x)   \rangle+\langle F_\alpha(x),D\varphi(x)  \rangle,\quad x\in H,
\end{equation}
where $D^2$ denotes the second Fr\'echet derivative.
Then, an easy application of It\^o's formula shows that the $L^1(H,\mu_\alpha)$-generator of $(P^\alpha_t)$
(given as before by $P^\alpha_tf(x)=\E [f(X_\alpha(t,x))]$) is given on
$\mathcal FC^2_b(\{e_n\})$ by $L_\alpha$. In particular,
$$
\int_{H}L_\alpha \varphi\;d\mu_\alpha=0,\quad\;\forall\;\varphi\in \mathcal FC^2_b(\{e_n\}).
$$
By a simple approximation argument and \eqref{e5.2} we get for $\alpha\in (0,\alpha_0)$ and
$$
\varphi_n(x):=\sum_{i=1}^nq_i^{-1}\langle x,e_i  \rangle^2,\quad x\in H,\;n\in \N,
$$
that also
\begin{equation}
\label{e5.21}
\int_{H}L_\alpha \varphi_n\;d\mu_\alpha=0.
\end{equation}
But for all $x\in H$, with $\tilde{F_\alpha}$ defined as $\tilde{F_0}$ in \eqref{e5.17}, we have
\begin{multline}
\label{e5.22}
L_\alpha \varphi_n(x)=-2\sum_{i=1}^n\frac{\lambda_i}{q_i}\langle x,e_i  \rangle^2+2\sum_{i=1}^nq_i^{-1}\langle
\tilde F_\alpha(x),e_i  \rangle\langle x,e_i  \rangle
\\+\sum_{i,j=1}^nq_i^{-1}\langle \sigma_n e_i,\sigma_ne_j  \rangle\\
\le -2\Theta(\pi_n x)+2\left(\sum_{i=1}^nq_i^{-1}\langle \tilde F_\alpha(x),e_i  \rangle^2\right)^{1/2}\;
\left(\sum_{i=1}^nq_i^{-1}\langle x,e_i  \rangle^2\right)^{1/2}
 \\+\sum_{i=1}^nq_i^{-1}\;|\sigma_ne_i|^2\\
 \le-2\Theta(\pi_nx)+c_1(1+|x|^{m+1}+\Theta^{1/2}(x)|x|)+\|\sigma\|^2\sum_{i=1}^\infty q_i^{-1},
\end{multline}
for some constant $c_1$ independent of $n$ and $\alpha.$ Here we used \eqref{e2.3} and \eqref{e5.18}. Now \eqref{e5.21}, \eqref{e5.2} and \eqref{e5.22} immediately imply that for some constant $\tilde c_1$
$$
\sup_{\alpha\in (0,\alpha_0)}\int_H\Theta(x)\mu_\alpha(dx)\le \sup_{\alpha\in (0,\alpha_0)}\tilde c_1\left(1+\int_{H}
|x|^{m+2}\mu_\alpha(dx)\right)+ \|\sigma\|^2\sum_{i=1}^\infty
q_i^{-1}<\infty.
$$
So, \eqref{e5.19} is proved, which by \eqref{e5.18}
implies \eqref{e5.4} and the proof is complete. $\Box$

 AKNOWLEDGEMENT. The second named author would like to thank UCSD, in particular,  his host Bruce Driver,
  for a very pleasant stay in La Jolla where a part of this work was done. The authors
  would like to thank Ouyang for his comments leading to a better constant
  involved in the Harnack inequality.

\end{document}